\documentclass[11pt,reqno]{amsart}

\usepackage[all,2cell]{xy}
\usepackage{amssymb}
\usepackage{mathrsfs}
\usepackage{cite}
\usepackage{amsfonts}
\newcommand{\rmv}[1]{}

\newcommand{\inv}{^{-1}}
\newcommand{\p}{\varphi}

\newcommand{\D}{{\mathscr D}} 

\newcommand{\R}{\mathrel{{\mathscr R}}} 
\newcommand{\eL}{\mathrel{{\mathscr L}}} 
\newcommand{\HH}{\mathrel{{\mathscr H}}} 

\newcommand{\malce}{\mathbin{\hbox{$\bigcirc$\rlap{\kern-8.3pt\raise0,50pt\hbox{$\mathtt{m}$}}}}}

\newcommand{\stab}[1]{\mathrm{St}(#1)}

\newcommand{\Thmname}{Theorem}
\newcommand{\Propname}{Proposition}
\newcommand{\Lemmaname}{Lemma}
\newcommand{\Definitionname}{Definition}
\newtheorem{Thm}{\Thmname}[section]
\newtheorem{Prop}[Thm]{\Propname}
\newtheorem{Lemma}[Thm]{\Lemmaname}
{\theoremstyle{definition}
\newtheorem{Def}[Thm]{\Definitionname}}
{\theoremstyle{remark}
}
\newtheorem{Cor}[Thm]{Corollary}

\newtheorem*{Lemma*}{Lemma}

\numberwithin{equation}{section}

\title{Semigroup Actions, Covering Spaces and Sch\"utzenberger Groups}

\author{Benjamin Steinberg}
\address{School of Mathematics and Statistics\\
Carleton University \\
1125 Colonel By Drive\\
Ottawa, Ontario  K1S 5B6 \\
Canada}
\thanks{The author was supported in part by NSERC and by a DFG Mercator guest professorship}
\email{bsteinbg@math.carleton.ca}
\date{May 29, 2009}

\dedicatory{In memory of Douglas Munn}

\keywords{Semigroups, Covering Spaces, Sch\"utzenberger groups}

\begin{document}

\begin{abstract}
We associate a $2$-complex to the following data: a presentation of a semigroup $S$ and a transitive action of $S$ on a set $V$ by partial transformations.  The automorphism group of the action acts properly discontinuously on this $2$-complex.  A sufficient condition is given for the $2$-complex to be simply connected.  As a consequence we obtain simple topological proofs of results on presentations of Sch\"utzenberger groups.  We also give a geometric proof that a finitely generated regular semigroup with finitely many idempotents has polynomial growth if and only if all its maximal subgroups are virtually nilpotent.
\end{abstract}

\maketitle

\section{Introduction}
The purpose of this article is to demonstrate that standard topological and geometric techniques from group theory can be used to study semigroups.  More precisely, given a semigroup presentation $\mathscr P=\langle X\mid R\rangle$ of a semigroup $S$ and a transitive right action of $S$ on a set $V$ by partial functions, we construct a connected 2-complex $K(V)$, called the action complex of $S$ on $V$.  If $G\leq \mathrm{Aut}_S(V)$, then the action of $G$ on $K(V)$ is without fixed-points and so there results a regular covering $\rho\colon K(V)\to G\backslash K(V)$.  We provide a sufficient condition to guarantee that $K(V)$ is simply connected, in which case $\pi_1(G\backslash K(V))\cong G$.  As a consequence, we obtain topological proofs of variations on well-known results concerning finite presentability of maximal subgroups and Sch\"utzenberger groups of semigroups~\cite{Ruskuc1,Ruskuc2}.  A topological approach along these lines for inverse semigroups was taken in~\cite{inversetop}.

The reader is referred to~\cite{CP,qtheor} for basic notions from semigroup theory such as Green's relations.  Let $S$ be a semigroup and denote by $S^1$ the result of adjoining an identity to $S$.  We recall that $s\R t$ (respectively, $s\eL t$) if $sS^1=tS^1$ (respectively, $S^1s= S^1t$).  The intersection of the relations $L$ and $\R$ is denoted $\HH$ and their join $\D$.  If $s\in S$, then $R_s$ denotes the $\R$-class of $s$.  Similar notation is used for the remainder of Green's relations.  If $R$ is an $\R$-class, then $\stab R=\{s\in S^1\mid sR\subseteq R\}$.  The quotient of $\stab R$ by its action on $R$ is a group $G(R)$ acting without fixed-points on $R$, known as the \emph{Sch\"utzenberger group} of $R$ (it depends only on the $\D$-class of $R$).  The orbits of $G(R)$ are precisely the $\HH$-classes of $R$; see~\cite{CP,qtheor}.  If $R$  contains an idempotent, then $G(R)$ is isomorphic to the maximal subgroup of $R$.  The \emph{right stabilizer} of $s\in S$ is $\stab s=\{t\in S\mid st=s\}$.

Let us say that a semigroup homomorphism is trivial if its image has cardinality at most $1$. Our first main theorem is the following result.

\begin{Thm}\label{main}
Let $S$ be a finitely presented semigroup and let $R$ be an $\R$-class of $S$ containing only finitely many $\HH$-classes.  Suppose that there exists $s\in R$ so that $\stab s$ admits no non-trivial homomorphism into a group. Then the Sch\"utzenberger group of $R$ is finitely presented.
\end{Thm}

Notice that if $S$ is left cancellative, then $\stab s$ is either trivial or empty and so Theorem~\ref{main} applies.  
Recall that if $s,t$ are elements of a semigroup $S$, one defines $s\eL^* t$ if, for all $x,y\in S^1$, one has $sx=sy$ if and only if $tx=ty$.  Trivially, $s\eL t$ implies $s\eL^* t$.  In particular, if $s$ is regular, then $s\eL^* e$ for some idempotent $e\in S$.  It is immediate that $s\eL^* t$ implies $\stab s=\stab t$.

\begin{Cor}
Let $S$ be a finitely presented semigroup and let $s\in S$ be such that the $\R$-class of $s$ has only finitely many $\HH$-classes and $\stab s=\stab e$ for some idempotent $e\in S$.   Then the Sch\"utzenberger group of $R_s$ is finitely presented.  The latter hypothesis applies in particular if $s$ is regular or, more generally, if the $\eL^*$-class of $s$ contains an idempotent.
\end{Cor}
\begin{proof}
An idempotent $e$ is a left zero of $\stab e$ and hence $\stab e$ has no non-trivial homomorphism into a group.  Theorem~\ref{main} now yields the desired conclusion.
\end{proof}

Theorem~\ref{main} encompasses the results of~\cite{Ruskuc1} and some of the results of~\cite{Ruskuc2}.  There is an example in~\cite{Ruskuc2} showing that having finitely many $\HH$-classes in an $\R$-class is not enough to ensure finite presentability of the Sch\"utzenberger group in general.  Of course, dual results apply to Sch\"utz\-en\-ber\-ger groups of $\eL$-classes (and we recall that the Sch\"utzenberger groups of $L_s$ and $R_s$ coincide). 

Recall that a semigroup is called \emph{right abundant} if each $\eL^*$-class contains an idempotent~\cite{Fountain}.  An immediate consequnce of the corollary is the following result.

\begin{Cor}
Let $S$ be a finitely presented right abundant semigroup (e.g.\ a regular semigroup).  Then if $R$ is an $\R$-class of $S$ containing only finitely many $\HH$-classes, the Sch\"utzenberger group of $R$ is finitely presented.
\end{Cor}
 
The other main result of this paper is that a finitely generated regular semigroup has polynomial growth if and only if all its maximal subgroups are virtually nilpotent.  This result is perhaps known to some --- certainly a version of this for linear groups can be found in~\cite{Okninski} --- but we believe this is the first geometric proof of this result.

\section{Action Complexes}
Let us first consider transitive actions of semigroups by partial functions.

\subsection{Automorphism groups of actions}
An \emph{action} of a semigroup $S$ on the right of a non-empty set $V$ by partial functions is a homomorphism $\p\colon S\to PT(V)$ where $PT(V)$ is the semigroup of partial transformations of $V$ (acting on the right).  As usual, we write $vs$ for $v\p(s)$.  The action is said to be \emph{transitive} if there are no non-empty, proper $S$-invariant subsets of $V$, or equivalently, $vS^1=V$ for all $v\in V$.  An \emph{automorphism} of the action is a bijective mapping $g\colon V\to V$ so that $(gv)s=g(vs)$ for all $s\in S$, where equality means that either both sides are defined and agree, or neither side is defined.  The group of automorphisms of the action is denoted $\mathrm{Aut}_S(V)$ 

The following lemma is immediate.

\begin{Lemma}\label{fixedpointfree}
Let $S$ act transitively on $V$ by partial maps.
Then the group $\mathrm{Aut}_S(V)$ acts on $V$ without fixed-points.
\end{Lemma}
\begin{proof}
Suppose that $g\in \mathrm{Aut}_S(V)$ fixes $v\in V$.  Let $w\in S\setminus \{v\}$.  Then $w=vs$ for some $s\in S$ and so $gw=g(vs)=(gv)s=vs=w$.  Thus $g=1$.  
\end{proof}

For instance, if $G$ is a group, then every subgroup of $G$ acts on the left of $G$ by automorphisms of the right regular representation of $G$.
The main example for us is the following.  Let $S$ be a semigroup and suppose that $R$ is an $\R$-class of $S$.  Then $S$ acts transitively on $R$ by partial functions via \[r\cdot s= \begin{cases} rs & rs\in R\\ \text{undefined} & \text{else.}\end{cases}\]
It is well known and easy to see that if $R$ is an $\R$-class of a semigroup $S$, then the Sch\"utzenberger group $G(R)$ acts by automorphisms on $R$~\cite{CP,qtheor}.  When $R$ is regular, then it is easy to show that $G(R)=\mathrm{Aut}_S(R)$.  Indeed, let $e\in R$ be an idempotent and suppose that $g$ is an automorphism with $ge=s$.  Then, for $r\in R$, one has $gr=ger=sr$.  In particular, $s\in \stab R$ and the actions on $R$ of the class of $s$ in $G(R)$ and of $g$ coincide.  

If $G\leq \mathrm{Aut}_S(V)$ is a group of automorphisms of the transitive action of $S$ on $V$, then $G\backslash V$ admits an induced transitive action of $S$ given by \[(Gv)s = \begin{cases} G(vs) & vs\ \text{is defined}\\ \text{undefined} & \text{otherwise.}\end{cases}\]  The natural map $\rho\colon V\to G\backslash V$ is a morphism of $S$-actions in the sense that $\rho(vs)=\rho(v)s$ for all $v\in V$ and $s\in S$, where again equality means that both sides are either undefined or are defined and agree.

\subsection{The action complex}
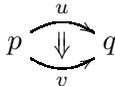
\begin{figure}[thbp]
\UseTwocells
\[\xymatrix{p \rtwocell^u_v&q}\]
\caption{Attaching the $2$-cells\label{2cells}}
\end{figure}

If $K$ is a $2$-complex (e.g.\ if $K$ is a graph), then the $i$-skeleton of $K$ will be denoted $K_i$.  Suppose now that $S$ acts on a set $V$ and that $\mathscr P=\langle X\mid R\rangle$ is a presentation of $S$.  If $w\in X^+$ is a word in the free semigroup on $X$, then $[w]_S$ will denote the corresponding element of $S$.  
The \emph{action complex} $K(\mathscr P,V)$ of $S$ on $V$ (written $K(V)$ if the presentation is understood) is the following $2$-complex.  The $0$-skeleton of $K(V)$ is $V$.  The edge set of $K(V)$ consists of all pairs $(v,x)\in V\times X$ such that $vx$ is defined.  We view the $1$-skeleton as a directed graph where $\iota(v,x)=v$ and $\tau(v,x)=vx$ determine the respective initial and terminal vertices of the edge $(v,x)$.  The set of $2$-cells consists of all pairs $(p,u=v)$ where $u=v\in R$, $p\in V$ and $q=p[u]_S=p[v]_S$ is defined.  The corresponding cell is defined as follows.  We take the unit disk and subdivide the northern semicircle from $(-1,0)$ to $(1,0)$ into $|u|$ equal length parts and the southern semicircle from $(-1,0)$ to $(1,0)$ into $|v|$ equal parts.  Then $(-1,0)$ is mapped to $p$, $(1,0)$ is mapped to $q$ and the northern (respectively, southern) semicircle is attached to the path labeled by $u$ (respectively $v$) from $p$.  See Figure~\ref{2cells}.
The $1$-skeleton $K(V)_1$ is termed the \emph{action graph} of $S$ on $V$ and depends only on the generating set.  The action graph for the action of a semigroup on one of its $\R$-classes is known in the literature as the \emph{Sch\"utzenberger graph} of the $\R$-class.

It is easy to see that transitivity of the action is equivalent to the directed graph $K(V)_1$ being strongly connected (i.e., any two points can be connected by a directed path).   A directed graph $T$ is said to be a \emph{directed tree rooted at a vertex $v$} if the underlying undirected graph of $T$ is a tree and the geodesic from $v$ to any vertex is a directed path.   The following lemma is a well-known application of Zorn's Lemma (cf.~\cite{fundgroupspaper,gst}).

\begin{Lemma}\label{directedspanningtree}
Let $\Gamma$ be a strongly connected directed graph and let $v$ be a vertex of $\Gamma$.  Then there is a directed spanning tree rooted at $v$ for $\Gamma$.
\end{Lemma}
\begin{proof}
Let $\mathscr C$ be the collection of all directed subtrees $T$ of $\Gamma$ rooted at $v$, ordered by inclusion.  It is non-empty since it contains $\{v\}$.  Next observe that if $\{T_{\alpha}\mid \alpha\in A\}$ is a chain in $\mathscr C$, then $\bigcup_{\alpha\in A} T_{\alpha}\in \mathscr C$.  Thus by Zorn's Lemma, $\mathscr C$ contains a maximal element $T$.  We assert that this is the desired directed spanning tree. Otherwise, since $\Gamma$ is strongly connected, there is a shortest directed path from $v$ to a vertex not in $T$. Let $e$ be the last edge of this path.  Then $\iota(e)\in T$ and $\tau(e)\notin T$.  One easily verifies that the graph $T'$ obtained by adjoining $e$ to $T$ is a directed tree rooted at $v$. 
Clearly $T'$ is a larger element of $\mathscr C$ than $T$.  This contradiction completes the proof.
\end{proof}

As a consequence, the fundamental group of a strongly connected directed graph is generated by directed loops~\cite{fundgroupspaper}. In what follows we do not distinguish notationally between paths and their homotopy classes.

\begin{Cor}\label{directedgenerate}
Let $\Gamma$ be a strongly connected directed graph and let $v\in \Gamma_0$ be a vertex.  Then $\pi_1(\Gamma,v)$ can be generated as a group by homotopy classes of directed loops.
\end{Cor}
\begin{proof}
Let $T$ be a directed spanning tree rooted at $v$.  For each $w\in \Gamma_0$, let $p_w$ be the geodesic path in $T$ from $v$ to $w$; by assumption it is directed.  Also choose a directed path $q_w$ from $w$ to $v$ (using strong connectivity).  Let $E$ be the set of positively oriented edges of $\Gamma$ not belonging to $T$.  We claim that $\pi_1(\Gamma,v)$ is generated by the homotopy classes of the directed paths of the form $p_{\iota(e)}eq_{\tau(e)}$ with $e\in E$ and $p_wq_w$ with $w\in V_0$.  Indeed, these paths represent elements of $\pi_1(\Gamma,v)$.  Moreover, if $e\notin E$, then \[p_{\iota(e)}ep_{\tau(e)}\inv = p_{\iota(e)}eq_{\tau(e)}(p_{\tau(e)}q_{\tau(e)})\inv\] (see Figure~\ref{positiveloops})
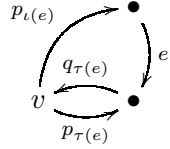
\begin{figure}[thbp]
\[\xymatrix{&\bullet\ar@/^/[d]^{e}\\ v\ar@/^1pc/[ur]^{p_{\iota(e)}}\ar@/_/[r]_{p_{\tau(e)}}&\bullet\ar@/_/[l]_{q_{\tau(e)}}}\]
\caption{Directed generators for $\pi_1(\Gamma,v)$\label{positiveloops}}
\end{figure}
and so each standard generator of $\pi_1(\Gamma,v)$ is in the subgroup generated by our collection of directed paths.
\end{proof}

Let us return to the situation of a semigroup $S=\langle X\mid R\rangle$ acting transitively by partial transformations on a set $V$.  Let $G\leq \mathrm{Aut}_S(V)$ be a subgroup of the automorphism group.  Then the action of $G$ on $V$ extends to a fixed-point free action of $G$ on $K(V)$ by defining $g(v,x)=(gv,x)$ on edges and $g(p,u=v) = (gp,u=v)$ on $2$-cells.  Moreover, notice that $G\backslash K(V)\cong K(G\backslash V)$, as is easily checked.  The quotient mapping $p\colon K(V)\to K(G\backslash V)$ is a regular covering map, as any fixed-point free automorphism group of a $2$-complex is properly discontinuous.  Consquently, 
\begin{equation}\label{regularcovering}
G\cong \pi_1(K(G\backslash V),Gv)/p_\ast(\pi_1(K(V),v).
\end{equation}
This leads to the following proposition.

\begin{Prop}
Let $S$ be a semigroup acting transitively on a set $V$ and let $G$ be a group of automorphisms of the action such that $G\backslash V$ is finite.  Then if $S$ is finitely generated, so is $G$.
\end{Prop}
\begin{proof}
The action graph $K(G\backslash V)$ with respect to a finite generating set is clearly finite.  Thus $G$ is finitely generated by \eqref{regularcovering}.
\end{proof}

By considering the action of the Sch\"utzenberger group on an $\R$-class, we obtain:

\begin{Cor}
Let $R$ be an $\R$-class of a finitely generated semigroup $S$ such that $R$ contains only finitely many $\HH$-classes.  Then the Sch\"utzenberger group $G(R)$ is finitely generated.
\end{Cor}

Next we provide a sufficient condition for $K(V)$ to be simply connected.  In this setting one obtains from \eqref{regularcovering} that $G\cong \pi_1(K(G\backslash V))$.

\begin{Def}[Stabilizer condition]
Let $S$ be a semigroup acting transitively on a set $V$.  We say that the action satisfies the \emph{stabilizer condition} if there is an element $v\in V$ so that $\stab v = \{s\in S\mid vs=v\}$ admits no non-trivial homomorphism into a group.
\end{Def}

Notice that if $S$ is a group, then the stabilizer condition is equivalent to the action being fixed-point free. Our aim is to show that the stabilizer condition implies that the action complex is simply connected.  We begin with a lemma.

\begin{Lemma}\label{positivepathsok}
Let $S=\langle X\mid R\rangle$ act on a set $V$ and let $w,w'\in X^+$ be such that $[w]_S=[w']_S$ and $p[w]_S=p[w']_S$ is defined for $p\in V$.  Then the respective paths labeled by $w$ and $w'$ from $p$ to $q=p[w]_S=p[w']_S$ are homotopic in $K(V)$.
\end{Lemma}
\begin{proof}
By an easy induction it suffices to consider the case $w'$ can be derived from $w$ in one step, i.e., there exists $u=v\in R$ and words $x,y\in X^+$ so that $w=xuy$ and $w'=xvy$.  Then $(p[x]_S,u=v)$ is a $2$-cell of $K(V)$ witnessing that the paths labeled by $u$ and $v$, respectively, from $p[x]_S$ to $p[xu]_S=p[xv]_S$ are homotopic; see Figure~\ref{homotopy}.
\begin{figure}[thbp]
\UseTwocells
\[\xymatrix{p\ar[r]^x&p[x]_S \rrtwocell^u_v&&p[xu]_S\ar[r]^y&q}\]
\caption{The homotopy of the paths labeled by $w$ and $w'$\label{homotopy}}
\end{figure}
Thus the respective paths labeled by $w$ and $w'$ from $[p]$ are homotopic, as required.
\end{proof}

We can now prove the main technical result of the paper.

\begin{Thm}\label{stabcond}
Let $S=\langle X\mid R\rangle$ be a semigroup equipped with a transitive action on a set $V$ satisfying the stabilizer condition.  Then $K(V)$ is simply connected.  Consequently, if $G\leq \mathrm{Aut_S(V)}$, then $G\cong  \pi_1(K(G\backslash V))$. 
\end{Thm}
\begin{proof}
The final statement is consequence of the first and \eqref{regularcovering}.  To prove the first statement, let $v\in V$ be as in the definition of the stabilizer condition. Let $T=\{w\in X^+\mid [w]_S\in \stab v\}$.  Then there results a homomorphism $\psi\colon T\to \pi_1(K(V),v)$ sending $w\in T$ to the homotopy class of the loop at $v$ labeled by $w$.  According to Corollary~\ref{directedgenerate}, the group $\pi_1(K(V),v)$ is generated by the classes of directed loops and hence $\psi(T)$ generates $\pi_1(K(V),v)$ as a group.  Lemma~\ref{positivepathsok} immediately yields  that $\psi$ factors through the projection \mbox{$T\to \stab v$}.  Thus $|\psi(T)|\leq 1$ (because every homomorphism from $\stab v$ into a group is trivial), and hence $\pi_1(K(V),v)$ is trivial, as required.
\end{proof}

Since the fundamental group of a finite $2$-complex is finitely presented and $K(G\backslash V)$ is finite whenever both $G\backslash V$ and $\langle X\mid R\rangle$ are finite, we obtain the following consequence of Theorem~\ref{stabcond}.

\begin{Cor}\label{finitepresentation}
Let $S$ be a finitely presented semigroup acting transitively on a set $V$ so that the stabilizer condition is satisfied.   Suppose $G$ is a group of automorphisms of the action with finitely many orbits.  Then $G$ is finitely presented.
\end{Cor}

A special case is when $H$ is a finite index subgroup of a finitely presented group $G$.  The action of $H$ on the left of $G$ is by automorphisms of the regular action of $G$ on the right of itself.  Since $H\backslash G$ is finite by assumption and a group is finitely presented as a group if and only if it is finitely presented as a semigroup, we recover the following classical result of Reidemeister.

\begin{Cor}[Reidemeister]
Finite index subgroups of finitely presented groups are finitely presented.
\end{Cor}

Theorem~\ref{main} is an immediate consequence of Corollary~\ref{finitepresentation} since the conditions of that theorem state that the action of the Sch\"utzenberger group on its $\R$-class has finitely many orbits and satisfies the stabilizer condition.

\section{Growth}
An undirected connected graph $\Gamma$ is \emph{locally finite} if the degree of each of its vertices is finite.  In this case we can turn the geometric realization of $\Gamma$ into a proper geodesic metric space by making each edge isometric to a unit interval.  The distance between two points is the length of the shortest path between them; this is called the \emph{path metric}.  Recall that a metric space is \emph{proper} if each closed ball is compact and is \emph{geodesic} if the distance between two points is the minimum length of a path between them.

Two metric spaces $(X,d)$ and $(X',d')$ are said to be \emph{quasi-isometric}~\cite{Bridson} if there is a function (called a \emph{quasi-isometry}) $f\colon X\to X'$ (not necessarily continuous) and $\lambda,\varepsilon,C>0$ so that:
\begin{itemize}
 \item $X'$ is contained in the $C$-neighborhood of $f(X)$;
\item For all $x,y\in X$, one has \[\frac{1}{\lambda}d(x,y)-\varepsilon\leq d'(f(x),f(y))\leq \lambda d(x,y)+\varepsilon.\]
\end{itemize}

Recall that if $G$ is a group generated as a group by a finite set $S$, then the \emph{word metric} on $G$ is defined by putting $d_G(g,h)$ to be the minimum length of a word representing $g\inv h$.  Equivalently, it is the metric induced by the path metric on the \emph{Cayley graph} of $G$ (i.e., the action graph of $G$ acting on the right of itself).  It is easy to verify that if $S,S'$ are two finite generating sets for $G$, then $(G,d_S)$ and $(G,d_{S'})$ are quasi-isometric~\cite{Bridson}.

A key tool in establishing quasi-isometries is the Milnor-\v{S}varc Lemma~\cite{Bridson}.  Recall that a group action of a group $G$ on a space $X$ is \emph{cocompact} if $G\backslash X$ is compact and is \emph{properly discontinuous} if, for each compact subset $K$ of $X$, there only finitely many elements $g\in G$ such that $gK\cap K\neq \emptyset$.  As was mentioned earlier, any fixed-point free action of a group on a $2$-complex is well known to be properly discontinuous.

\begin{Lemma}[Milnor-\v{S}varc]\label{milnor}
Suppose that $G$ is a group acting properly discontinuously and cocompactly by isometries on a proper geodesic metric space $X$.  Then $G$ is finitely generated and, for any basepoint $x_0$ of $X$, the map $g\mapsto gx_0$ is a quasi-isometry of $G$ and $X$.
\end{Lemma}

As a consequence, let us show that a group of automorphisms of a transitive action of a finitely generated semigroup with finitely many orbits is quasi-isometric to the action graph (which is a proper geodesic metric space).  

\begin{Prop}\label{localfiniteness}
Let $S$ be a semigroup generated by a finite set $X$ acting transitively on a set $V$.  Suppose that $V$ admits a group $G$ of automorphisms with $G\setminus V$ finite.  Then the action graph $K(V)_1$ is locally finite, $G$ is finitely generated and, for any vertex $v_0\in V$, there is a quasi-isometry $G\to K(V)_1$ given by $g\mapsto gv_0$.
\end{Prop}
\begin{proof}
The projection $\p\colon K(V)_1\to K(G\setminus V)_1$ is a covering and so preserves the degree of a vertex.  But $K(G\setminus V)_1$ is a finite graph and so each vertex has finite degree.  Thus $K(V)_1$ is locally finite and hence a proper geodesic metric space with respect to the path metric.  The action of $G$ on $K(V)_1$ is by graph automorphisms and hence by isometries of the path metric.  The quotient $G\backslash K(V)_1=K(G\backslash V)_1$ is a finite graph, and hence compact, so that action of $G$ is cocompact. Also, since the action is fixed-point free, it is properly discontinuous.  The remainder of the proposition is now immediate from the Milnor-\v{S}varc Lemma. 
\end{proof}

Recall that there is a preorder on non-decreasing functions from the natural numbers to the non-negative reals given by $f\preceq g$ if there is a constant $k$ so that $f(n)\leq kg(kn+k)+k$ for all $n\in N$.  One write $f\approx g$ if $f\preceq g$ and $g\preceq f$.  In this case we say that $f$ and $g$ have the same \emph{order of growth}.  We say that $f$ has \emph{polynomial growth of degree at most $d$ (of degree $d$)} if $f\preceq n^d$ ($f\approx n^d$).  In general, $f$ is said to have \emph{polynomial growth} if $f$ has polynomial growth of degree at most $d$ for some $d$.

The \emph{growth function} of semigroup $S$ with respect to a finite generating set $X$ is the mapping $g_{S,X}\colon \mathbb N\to \mathbb R$ given by $g_{S,X}(n) = |S_n|$ where \[S_n=\{s\in S\mid s=[w]_S\ \text{with}\ |w|\leq n\}.\]  It is well known and easy to prove that if $X'$ is another finite generating set for $X$, then $g_{S,X}\approx g_{S,X'}$ and so in particular the order of growth of $S$ is well defined.  Also observe that if $T$ is a finitely generated subsemigroup of $S$, then the order of growth of $T$ is bounded above by the order of growth of $S$ (just adding a generating set of $T$ to the generating set of $S$).   A semigroup is said to have \emph{polynomial growth} if its growth function has polynomial growth.    A celebrated result of Gromov~\cite{Gromov} says that a finitely generated group $G$ has polynomial growth if and only if it is \emph{virtually nilpotent}, that is, has a nilpotent subgroup of finite index.  The extension of this result to cancellative semigroups was obtained by Grigorchuk~\cite{Grigpoly}.

If $\Gamma$ is a locally finite undirected graph and $v$ is a vertex $\Gamma$, then the growth function of $\Gamma$ at $v$ is defined by $g_{\Gamma,v}(n) = |B_{\Gamma}(v,n)\cap \Gamma_0|$ where $B_{\Gamma}(v,n)$ is the closed ball of radius $n$ around $v$.  It is well known (cf.~\cite{Bridson}) that if $f\colon \Gamma\to \Gamma'$ is a quasi-isometry, then $g_{\Gamma,v}\approx g_{\Gamma',f(v)}$.  In particular, the growth rate of a group is the same as that of its Cayley graph 

Suppose now that $\Gamma$ is a directed graph with finite outdegree at each vertex.  There is a very natural growth function associated to a vertex $v$ of $\Gamma$, called the \emph{directed growth function}, given by putting $\vec{g}_{\Gamma,v}(n)$ equal to the number of elements reachable from $v$ by a directed path of length at most $n$.  For instance, the growth rate of a finitely generated monoid is the same as the directed growth rate of its Cayley graph.  It turns out that for action graphs of the sort we have been discussing, these growth functions are equivalent.

\begin{Prop}\label{reverseedges}
Let $S$ be a semigroup generated by a finite set $X$ acting transitively on a set $V$.  Suppose that $V$ admits a group $G$ of automorphisms with $G\setminus V$ finite.  Let $v$ be an element of $V$.  Then $g_{K(V)_1,v}\approx \vec g_{K(V)_1,v}\preceq g_{S,X}$ and $g_{K(V)_1,v}$ has the same order of growth as $G$.   Consequently, if $S$ has polynomial growth, then $G$ is a finitely generated virtually nilpotent group.
\end{Prop}
\begin{proof}
Clearly $\vec g_{K(V)_1,v}\preceq g_{K(V)_1,v}$.  For the converse, we observe that there are only finitely many orbits of edges in $K(V)$ (equal to the number of edges of $K(G\backslash V)$).  Fix edges $e_1,\ldots,e_n$ representing the distinct orbits.  Then for each $i$, we can find a directed path $p_i$ from $\tau(e_i)$ to $\iota(e_i)$ by strong connectivity.  Let $k$ be bigger than the length of any of the $p_i$.  Then, for any edge $e$ of $K(V)_1$, there is a directed path of length at most $k$ from $\tau(e)$ to $\iota(e)$.  Hence each vertex of $B_{K(V)_1}(v,n)$ can be reached by a directed path of length at most $kn$.   Indeed, if $p$ is an undirected path of length at most $n$ from $v$ to $w$, then we can replace every negatively traversed edge by a directed path of length at most $k$.  Thus $g_{K(V)_1,v}(n)\leq \vec g_{K(V)_1,v}(kn)$ and so $g_{K(V)_1,v}\approx \vec g_{K(V)_1,v}$.  

There is a surjective partial map from $S_n\cup \{1\}$ to the set of elements reachable from $v$ by a directed path of length at most $n$ given by $s\mapsto vs$ if $vs$ is defined, and otherwise the map is undefined on $s$.  It then follows that $\vec g_{K(V)_1,v}\preceq g_{S,X}$.

The Milnor-\v{S}varc Lemma immediately implies that the order of growth of $G$ is the same as the order of growth of $g_{K(V)_1,v}$.  The final statement is a consequence of Gromov's theorem on groups of polynomial growth.
\end{proof}

Applying the above result to Sch\"utzenberger groups we obtain:
\begin{Cor}\label{Schutzgraph}
Let $S$ be a finitely generated semigroup of polynomial growth at most $d$ and let $R$ be an $\R$-class of $S$ containing only finitely many $\HH$-classes.  Then the Sch\"utzenberger group of $R$ is a finitely generated virtually nilpotent group of polynomial growth of degree at most $d$.
\end{Cor}

Of course, the above corollary is trivial for maximal subgroups since they are finitely generated subsemigroups of $S$; but the Sch\"utzenberger groups of non-regular $\R$-classes are not subsemigroups of $S$ so one must give some argument.  
It turns out that the converse of Corollary~\ref{Schutzgraph} is true for regular semigroups with finitely many idempotents.  We believe this result to be new, although a similar result for linear semigroups can be found in~\cite{Okninski}.  In any event, the point is to provide a simple geometric argument.

\begin{Thm}
Let $S$ be a finitely generated regular semigroup with finitely many idempotents.  Then $S$ has polynomial growth of degree at most $d$ if and only if each of its maximal subgroups is virtually nilpotent of polynomial growth of degree at most $d$.
\end{Thm}
\begin{proof}
A regular semigroup has finitely many idempotents if and only if it has finitely many $\eL$- and $\R$-classes.  In particular, each maximal subgroup acts on the Sch\"utzenberger graph of its $\R$-class with finitely many orbits and so the above results apply (and in particular, each maximal subgroup is finitely generated).  Hence if $S$ has polynomial growth of degree at most $d$, then each maximal subgroup of $S$ is virtually nilpotent with polynomial growth of degree at most $d$ (say by an application of Corollary~\ref{Schutzgraph}).

Suppose conversely that each maximal subgroup of $S$ is virtually nilpotent of polynomial growth of degree at most $d$.     Let $e_1,\ldots, e_m$ be a complete set of (distinct) idempotent representatives of the $\R$-classes of $S$ and denote by $\Gamma_i$ the Sch\"utzenberger graph of the $\R$-class of $e_i$.  Proposition~\ref{reverseedges} implies that each growth function $g_{\Gamma_i,e_i}$ has polynomial growth of degree at most $d$.  Let $s\in S_n$ with $s=[w]_S$ where $|w|\leq n$.  By regularity, $s\R e_i$ for a unique $e_i$.  As $e_is=s$, it follows that $w$ labels a path from $e_i$ to $s$ in $\Gamma_i$ and hence $s\in B_{\Gamma_i}(e_i,n)$.  Thus we have defined an injective function $S_n\to \biguplus_{i=1}^m B_{\Gamma_i}(e_i,n)$ sending $s$ to the vertex $s$ of the Sch\"utzenberger graph of its $\R$-class and so \[g_{S,X}(n)\leq \sum_{i=1}^mg_{\Gamma_i,e_i}(n).\]  We conclude $g_{S,X}\preceq x^d$, as required.
\end{proof}

The assumption on finiteness of the number of idempotents is needed in this theorem.  For instance, in~\cite[Example 5.3]{okninskiadv} a $3$-generated completely regular semigroup $S$ is constructed consisting of an infinite cyclic group of units and a completely $0$-simple $0$-minimal ideal with maximal subgroup free nilpotent of class $2$ such that $S$ contains a finitely generated subsemigroup of intermediate growth (i.e., of growth faster than any polynomial but subexponential).

\subsection*{Acknowledgments}
When I was a graduate student, John Rhodes said to me that group mapping is a covering space of right letter mapping.  This remark inspired~\cite{inversetop} and this paper.

\def\malce{\mathbin{\hbox{$\bigcirc$\rlap{\kern-7.75pt\raise0,50pt\hbox{${\tt
  m}$}}}}}\def\cprime{$'$} \def\cprime{$'$} \def\cprime{$'$} \def\cprime{$'$}
  \def\cprime{$'$} \def\cprime{$'$} \def\cprime{$'$} \def\cprime{$'$}
  \def\cprime{$'$}


\begin{thebibliography}{10}

\bibitem{Bridson}
M.~R. Bridson and A.~Haefliger.
\newblock {\em Metric spaces of non-positive curvature}, volume 319 of {\em
  Grundlehren der Mathematischen Wissenschaften [Fundamental Principles of
  Mathematical Sciences]}.
\newblock Springer-Verlag, Berlin, 1999.

\bibitem{okninskiadv}
F.~Ced{\'o} and J.~Okni{\'n}ski.
\newblock Semigroups of matrices of intermediate growth.
\newblock {\em Adv. Math.}, 212(2):669--691, 2007.

\bibitem{CP}
A.~H. Clifford and G.~B. Preston.
\newblock {\em The algebraic theory of semigroups. {V}ol. {I}}.
\newblock Mathematical Surveys, No. 7. American Mathematical Society,
  Providence, R.I., 1961.

\bibitem{Fountain}
J.~Fountain.
\newblock Abundant semigroups.
\newblock {\em Proc. London Math. Soc. (3)}, 44(1):103--129, 1982.

\bibitem{Grigpoly}
R.~I. Grigorchuk.
\newblock Semigroups with cancellations of degree growth.
\newblock {\em Mat. Zametki}, 43(3):305--319, 428, 1988.

\bibitem{Gromov}
M.~Gromov.
\newblock Groups of polynomial growth and expanding maps.
\newblock {\em Inst. Hautes \'Etudes Sci. Publ. Math.}, (53):53--73, 1981.

\bibitem{gst}
J.~McCammond, J.~Rhodes, and B.~Steinberg.
\newblock Geometric semigroup theory.
\newblock In preparation.

\bibitem{Okninski}
J.~Okni{\'n}ski.
\newblock {\em Semigroups of matrices}, volume~6 of {\em Series in Algebra}.
\newblock World Scientific Publishing Co. Inc., River Edge, NJ, 1998.

\bibitem{qtheor}
J.~Rhodes and B.~Steinberg.
\newblock {\em The {$q$}-theory of finite semigroups}.
\newblock Springer Monographs in Mathematics. Springer, New York, 2009.

\bibitem{Ruskuc1}
N.~Ru{\v{s}}kuc.
\newblock Presentations for subgroups of monoids.
\newblock {\em J. Algebra}, 220(1):365--380, 1999.

\bibitem{Ruskuc2}
N.~Ru{\v{s}}kuc.
\newblock On finite presentability of monoids and their {S}ch\"utzenberger
  groups.
\newblock {\em Pacific J. Math.}, 195(2):487--509, 2000.

\bibitem{fundgroupspaper}
B.~Steinberg.
\newblock Fundamental groups, inverse {S}ch\"utzenberger automata, and monoid
  presentations.
\newblock {\em Comm. Algebra}, 28(11):5235--5253, 2000.

\bibitem{inversetop}
B.~Steinberg.
\newblock A topological approach to inverse and regular semigroups.
\newblock {\em Pacific J. Math.}, 208(2):367--396, 2003.

\end{thebibliography}
\end{document}